\documentclass[12pt]{article}

\usepackage{amsmath,amssymb}
\usepackage[mathscr]{eucal}

\newcommand{\pl}{\partial}
\newcommand{\ol}{\overline}
\newcommand{\plr}{\partial r}
\newcommand{\tv}{\tilde{v}}
\newcommand{\tu}{\tilde{u}}
\newcommand{\tw}{\tilde{w}}
\newcommand{\Dtph}{\Delta_{\Omega}}
\newcommand{\DD}{\mathcal{D}}
\newcommand{\SP}{\mathcal{S}}

\newcommand{\Ww}{\mathcal{W}}
\newcommand{\RE}{\mathrm{Re}}
\newcommand{\RR}{\mathbb{R}}
\newcommand{\CC}{\mathbb{C}}
\newcommand{\Sph}{\mathbb{S}}
\newcommand{\OO}{\mathcal{O}}

\begin{document}
\begin{center}
{\Large\bf Extensions of the quadratic form\\[3mm]
           of the transverse Laplace operator}\\[1cm]

{T.~A.~Bolokhov}

\vspace{0.6cm}
{\it St.Petersburg Department of V.\,A.\,Steklov Mathematical Institute\\ 
         Russian Academy of Sciences\\
         27 Fontanka, St.\,Petersburg, Russia 191023}

\end{center}
\vspace{0.6cm}


\begin{abstract}
    We review the quadratic form of the Laplace operator in 3 dimensions
    written in spherical coordinates and acting on the transverse components of 
    vector functions.
    Operators, acting on the parametrizing functions of one of the transverse
    components with angular momentum 1 and 2,
    appear to be fourth order symmetric operators
    with deficiency indices (1,1).
    We develop self-adjoint extensions of these operators
    and propose correspondent extensions for the initial quadratic form.
    The relevant scalar product for the angular momentum 2
    differs from the original product in the space of vector functions,
    but nevertheless it is still local in radial variable.
    Eigenfunctions of the operator extensions in question can be
    treated as stable soliton-like solutions of the corresponding
    dynamical system with the quadratic form being a functional
    of the potential energy.
\end{abstract}

\section*{Introduction}
    The quadratic form of the Laplace operator acting on the transverse
    components of a vector field
$ \vec{f}(\vec{x}) $
\begin{equation}
\label{TVL}
    Q_{0}(\vec{f}) =
    (\vec{f}, \Delta \vec{f}) =
    - \sum_{j,k=1}^{3} \int_{\RR^{3}}
	\ol{f^{k}} \frac{\pl^{2}}{\pl x_{j}^{2}} f^{k} d^{3}x,
    \quad \vec{\pl}\cdot\vec{f} = \sum_{k} \frac{\pl f^{k}}{\pl x_{k}} = 0
\end{equation}
    plays the role of the energy functional in the Coulomb gauge
    electrodynamics.
    The Friedrichs-Stone theorem
\cite{FS}
    states, that each closed semi-bounded
    quadratic form in a Hilbert space corresponds to a self-adjoint operator.
    That is, in order to define an operator for the theory (and altogether
    the spectral decomposition and the dynamics) one has to define a scalar
    product and a dense domain on which the form is acting.
    This is not always straightforward and usually requires taking into
    account
    additional information about physical properties of the model.
    Moreover, one can include a part of the scalar product into the operator
    expression (and vice versa) and in this way possibly find
    a different dynamics.

    Even when the scalar product and the domain are both fixed there is
    still a possibility that the form in question can be extended to a
    wider domain. This happens when the formal operator
    of the quadratic form, defined on the initial domain of the problem,
    is symmetric but not self-adjoint. An example of such a model is
    one with the quadratic form of the scalar Laplace operator in two- or
    three-dimensional space
\begin{equation}
\label{SL}
    Q_{0}(f) = (f, \Delta f) = -\sum_{j=1}^{3} \int_{\RR^{3}} \ol{f(\vec{x})}
	\frac{\pl^{2}}{\pl x_{j}^{2}} f(\vec{x}) d^{3}x .
\end{equation}
    Even in the ``plain'' scalar product
\begin{equation*}
    (f,g) = \int_{\RR^{3}} \ol{f(\vec{x})} g(\vec{x}) d^{3}x ,
\end{equation*}
    written in spherical coordinates, the operator acting on
$ s $-subspace (the subspace corresponding to spherically symmetric
    functions with angular momentum zero) turns out to be symmetric
    with deficiency indices
$ (1,1) $.
    And the corresponding quadratic form
$ Q_{0}(f) $
    can be extended to the
    subset of functions diverging at zero as
$ |x|^{-1} $.

    Quadratic form
(\ref{TVL})
    is somewhat different.
    Parametrization of the transverse component of the vector field with
    two sets of functions reveals, that the quadratic forms of one of them are
    represented by differential operators of the order 4. In the subspaces
    corresponding to the angular momentum 
$ l=1 $ and
$ l=2 $
    these operators are symmetric with deficiency indices 
$ (1,1) $.
    We shall construct self-adjoint extensions of these operators and of the
    corresponding quadratic forms. From the start these extensions are
    closed with repect to the ``plain'' product in the space of parametrizing
    functions, which naturally differs from the product induced by
\begin{equation}
\label{fgprod}
    (\vec{f}, \vec{g}) =
	\int_{\RR^{3}} \ol{\vec{f}(\vec{x})} \cdot \vec{g}(\vec{x}) d^{3}x .
\end{equation}
    Nevertheless in the case of the angular momentum
$ l=1 $
    domains of the extensions in questions
    are also finite (and expected to be closed) with respect to the latter
    product.

    The article is composed as follows. Section~\ref{sscalar} provides
    some classical background on extensions of quadratic form of the
    scalar Laplace operator.
    Sections \ref{svector}, \ref{strav} and \ref{soper}
    introduce the basis for vector
    spherical harmonics, the parametrization of the transverse subspace
    and the expressions for quadratic forms.
    Sections \ref{sdef}, \ref{ssymm}, \ref{sdiscrete} and \ref{scont}
    describe self-adjoint extensions of the fourth-order
    differential operators and their spectral decompositions. 
    The normalization of the eigenfunctions and the completeness relations
    are given in the Section \ref{snorm}. Kernels of the resolvents and of
    the reverse operators are in the Section \ref{resolvent}.
    Section \ref{squadr} provides description of the extended
    quadratic form of the 3-dimensional Laplace operator.
$ l=2 $.

\section{Scalar Laplace operator}
\label{sscalar}
    The scalar Laplace operator acts on the space of twice-differentiable 
    functions of three variables as follows
\begin{equation*}
    \Delta: \quad f(\vec{x}) \to \Delta f(\vec{x}) = -\sum_{j=1}^{3}
	\frac{\pl^{2} f}{\pl x_{j}^{2}}\,.
\end{equation*}
    Switching to the spherical coordinates
\begin{gather*}
    \vec{x} = \vec{x}(r,\theta,\varphi)
    = \begin{pmatrix} r\cos\theta \cos\varphi\\
	r\cos\theta \sin\varphi\\
	r\sin\theta
	\end{pmatrix}, \\
    0 \leq r,\quad 0\leq\theta\leq\pi,\quad 0 \leq\varphi < 2\pi
\end{gather*}
    brings the operator to the form
\begin{equation*}
    \Delta f(\vec{x}(r,\theta,\varphi)) =
	-\frac{1}{r^{2}}\frac{\pl}{\plr} r^{2} \frac{\pl}{\plr} f
	-\frac{1}{r^{2}\sin\theta} \bigl(
    \frac{\pl}{\pl\theta}\sin\theta \frac{\pl}{\pl\theta}
    + \frac{1}{\sin\theta}\frac{\pl^{2}}{\pl \varphi^{2}}\bigr) f .
\end{equation*}
    The standard separation of variables
\begin{equation}
\label{scalarSep}
    f(r,\theta,\varphi)
	= \sum_{0\leq |m|\leq l} f_{lm}(r) Y_{lm}(\theta,\varphi)
\end{equation}
    allows one to transform the action of the operator to
\begin{equation}
\label{Lsc}
    \Delta f = \sum_{0\leq |m|\leq l} \Bigl(
	- \frac{1}{r^{2}}\frac{\pl}{\plr} r^{2} \frac{\pl}{\plr}
	+ \frac{l(l+1)}{r^{2}} \Bigr) f_{lm}(r) Y_{lm}(\theta,\varphi) ,
\end{equation}
    where
$ Y_{lm}(\theta,\varphi) $ are the spherical harmonics satisfying
\begin{gather*}
    \Delta_{\Omega} Y_{lm} = 
	- \frac{1}{\sin\theta} \Bigl(
    \frac{\pl}{\pl\theta}\sin\theta \frac{\pl}{\pl\theta}
    + \frac{1}{\sin\theta}\frac{\pl^{2}}{\pl \varphi^{2}}\Bigr) Y_{lm}
    = l(l+1) Y_{lm}\,, \\
    \int_{\Sph^{2}} \ol{Y_{l'm'}(\Omega)} Y_{lm}(\Omega) d\Omega
	= \delta_{ll'} \delta_{mm'} .
\end{gather*}
    Since in what follows the equations will not depend on the
    parametrization of a point on sphere 
$ \Sph^{2} $
    by spherical angles
$ \theta $, $ \varphi $,
    we have replaced the latter with a generic spherical variable
$ \Omega $.
    For our purposes, it is only the completeness of the set
$ Y_{lm}(\Omega) $
    that matters, 
    {\it i.\,e.} we require that every smooth function be
    uniquely represented as sum
(\ref{scalarSep}).
    Thus the action of the Laplace operator reduces to the action of the
    radial operators
\begin{equation}
\label{TauL}
    \tau_{l} : \quad f_{l}(r) \to
	-\frac{1}{r^{2}}\frac{d}{dr}r^{2}\frac{d}{dr}f_{l}(r)
	    + \frac{l(l+1)}{r^{2}} f_{l}(r)
\end{equation}
    on the space of differentiable functions defined on the positive
    half-axis, endowed with the scalar product
\begin{equation*}
    (f,g)_{\RR^{3}} = \int_{0}^{\infty} \ol{f(r)} g(r)\, r^{2} dr .
\end{equation*}
    By virtue of functional substitution
$ f_{l} = \frac{u_{l}}{r} $
    the operators
$ \tau_{l} $
    transform into
\begin{equation}
\label{Tl}
    T_{l} : \quad u_{l} \to - \frac{d^{2}}{dr^{2}} u_{l}
	+ \frac{l(l+1)}{r^{2}} u_{l} ,
\end{equation}
    while the scalar product turns into a simpler expression
\begin{equation}
\label{plainprod}
    (v,u) = \int_{0}^{\infty} \ol{v(r)} u(r) dr .
\end{equation}

\subsection{Self-adjointness of the radial operators}
\label{saradial}
    Operators
$ T_{l} $,
    defined on the subspace of twice-differentiable functions vanishing
    along with their derivatives at zero
\begin{equation*}
    \Ww_{0}^{2} = \{u:\quad (u,u)<\infty,\, (u'',u'')<\infty,
	\, u(0)=u'(0)=0 \} ,
\end{equation*}
    are essentially self-adjoint for
$ l\geq 1 $.
    Meanwhile the operator
$ T_{0} $ (from now on we shall denote it as
$ T $)
    has the deficiency indices
$ (1,1) $.
    Indeed, if
$ T $
    is defined on
$ \Ww_{0}^{2} $,
    then the functions
\begin{equation*}
    g_{\pm} = \exp\{e^{\mp i\frac{3\pi}{4}} \rho r\}
\end{equation*}
    generate the kernels of adjoint operators
$ (T\mp i\rho^{2})^{*} $,
    since they are the only square integrable
    solutions of the equations
\begin{equation}
\label{Tadj}
    \frac{d^{2}g_{+}}{dr^{2}} = i\rho^{2} g_{+} \quad \text{ and } \quad 
    \frac{d^{2}g_{-}}{dr^{2}} = - i\rho^{2} g_{-} .
\end{equation}
    Here we have introduced a dimensional parameter
$ \rho $: 
$ [\rho] = [r]^{-1} $
    for the reason that
$ T $ is a dimensional operator, {\it i.\,e.} its spectrum is comprised by
    dimensional values.
    From the point of view of physics, one cannot equate
    the left- and right-hand sides of
(\ref{Tadj})
    without introducing a dimensional coefficient of
$ \rho^{2} $.
    The choice of value of
$ \rho $
    in some sense is equivalent to the choice of the renormalization
    point in physics, being a quantity via which 
    a dimensional parameter of the theory can be expressed.

\subsection{Extensions of symmetric operator
$ T $}
    Exploiting the symmetric property of
$ T $,
    the ranges
\begin{equation*}
    Ran_{\pm} = \{g = (T\mp i\rho^{2})u:\enskip u \in \Ww_{0}^{2} \}
\end{equation*}
    can be isometrically connected by the operator
$ U $ (called the Cayley transform of 
$ T $)
    which acts as follows
\begin{equation*}
    U : \quad t= (T+i\rho^{2})u \to Ut = (T-i\rho^{2})u .
\end{equation*}
    Using the fact that the linear spans of vectors
$ g_{\pm} $
    represent the only
    orthogonal complements of the corresponding ranges
$ Ran_{\pm} $
    (by the definition of an adjoint operator),
    one may extend the partial isometry
$ U $
    to an unitary operator
$ U_{a} $
    by fixing it on
$ g_{+} $
    in the following way
\begin{equation*}
    U_{a} e^{ia} g_{+} = e^{-ia} g_{-} , \quad 0\leq a < \pi ,
\end{equation*}
    where
$ e^{2ia} $ is an unitary parameter.
    The operator
$ U_{a} $
    is a Cayley transform of some self-adjoint extension
$ T_{a} $
    of the symmetric operator
$ T $.
    This extension is defined on the subspace
$ \Ww_{a}^{2} $,
    which includes the function
$ h^{a} $
    obeying
\begin{align*}
    (T_{a}-i\rho^{2}) h^{a} &= e^{ia} g_{+} ,\\
    (T_{a}+i\rho^{2}) h^{a} &= e^{-ia} g_{-} ,
\end{align*}
    and the following relation holds for that subspace
\begin{equation*}
    \Ww_{a}^{2} = \Ww_{0}^{2} \dotplus \{\alpha h^{a}, \alpha\in\CC \} .
\end{equation*}
    It is not hard to see that
$ h^{a} $
    can be constructed as a linear combination of
$ g_{+} $ and
$ g_{-} $:
\begin{equation*}
    h^{a} = \frac{e^{ia}}{2i\rho^{2}} g_{+} - \frac{e^{-ia}}{2i\rho^{2}}g_{-}
	= \frac{1}{2i\rho^{2}} \bigl(\exp\{ia + e^{-i\frac{3\pi}{4}}\rho r\}
    - \exp\{-ia + e^{i\frac{3\pi}{4}}\rho r\} \bigr) .
\end{equation*}

\subsection{Boundary conditions}
    The operator of the second derivative is symmetric
    on any linear subspace on which the boundary terms
\begin{equation*}
    \int_{0}^{\infty} \ol{v''} u dr -
    \int_{0}^{\infty} \ol{v} u'' dr
	= (\ol{v'}u - \ol{v}u')|_{0}^{\infty} = \ol{v}(0)u'(0)
	    - \ol{v'}(0) u(0) 
\end{equation*}
    vanish (here we have eliminated the limits at infinity using the square
    integrability of
$ u $, $ v $).
    The element
$ h^{a} $
    constructed above relates
    the value of function
$ u $ from
$ \Ww_{a}^{2} $
    at zero
    with the value of its derivative in such a way that the operator
$ T_{a} $
    is symmetric:
\begin{equation*}
    \ol{h_{a}(0)} u'(0) - \ol{h'_{a}(0)} u(0) = 0 , 
\end{equation*}
    and hence
\begin{equation*}
    \Ww_{a}^{2} = \{u: (u,u)<\infty, \: (u'',u'')<\infty,\:
	\rho\sin(a-\frac{\pi}{4})\, u(0) +\sin a\, u'(0) = 0 \} .
\end{equation*}
    Essential self-adjointness of
$ T_{a} $
    on
$ \Ww_{a}^{2} $
    (deficiency indices of
$ (0,0) $) follows from the construction.

\subsection{Eigenvectors}
    Depending on the value of parameter
$ a $,
    the operator
$ T_{a} $
    might have eigenvalues (discrete spectrum).
    Consider the exponential
\begin{equation*}
    v_{\kappa}(r) = e^{-\kappa r}, \quad \kappa > 0 ,
\end{equation*}
    which obeys the equation
\begin{equation*}
    - \frac{d^{2}v_{\kappa}}{dr^{2}} = - \kappa^{2} v_{\kappa} .
\end{equation*}
    This exponential belongs to 
$ \Ww_{a}^{2} $
    if the following condition holds
\begin{equation*}
    \rho \sin(a-\frac{\pi}{4})\, v_{\kappa}(0) = - \sin a \, v_{\kappa}{}'(0).
\end{equation*}
    This condition relates
$ \kappa $
    with parameters
$ a $ and
$ \rho $ as follows
\begin{equation*}
    \kappa = \rho \frac{\sin(a-\pi/4)}{\sin a} .
\end{equation*}
    As
$ a $
    is varied from
    0 to
$ \pi $,
    the RHS runs through the values of the real axis once.
    It is clear, that in order for 
$ \kappa $
    to be positive, the numerator above should be positive, {\it i.e.}
\begin{equation*}
    \frac{\pi}{4} < a < \pi .
\end{equation*}
    There is no discrete spectrum of
$ T_{a} $
    for other values of 
$ a $.
    The continuous spectrum of
$ T_{a} $
    has multiplicity one and
    occupies negative half-axis, we shall not stop on it here, thorough
    material about this operator and the theory of self-adjoint extensions
    can be found in
\cite{RS}.

\subsection{The quadratic form}
    The Friedrichs-Stone theorem
\cite{FS}
    uniquely relates each semi-bounded self-adjoined operator
$ T_{a} $
    to a closed semi-bounded quadratic form
$ Q_{a} $,
    via a natural rule
\begin{equation*}
    Q_{a}(u) = (u,T_{a}u), \quad u \in \Ww_{a}^{2} .
\end{equation*}
    This quadratic form is initially defined on the domain
$ \Ww_{a}^{2} $,
    and then continuously extended onto some closed domain
$ \Ww_{(a)}^{1} $.
    The Friedrichs-Stone theorem also ensures that among 
    self-adjoint extensions of the symmetric semi-bounded operator
$ T $
    (with deficiency indices
$ (1,1) $ in this case)
    there is a distinguished (maximal) extension whose quadratic form
    is the closure of the quadratic form of
$ T $:
\begin{equation*}
    Q_{0}(u) = (u,Tu) , \quad u \in \Ww_{0}^{2} .
\end{equation*}
    In the case of the second derivative operator this distinguished
    (Friedrichs) extension corresponds to
$ a=0 $,
$ \kappa=-\infty $ and the domain of the corresponding form
$ Q_{0} $
    consists of functions vanishing at zero:
\begin{equation*}
    \Ww_{0}^{1} = \{u: \enskip (u,u)<\infty,\, (u',u')<\infty,\, u(0)=0 \} .
\end{equation*}
    All other self-adjoined extensions correspond to forms
    defined on one single domain
\begin{equation*}
    \Ww_{(a)}^{1} = \Ww_{1}^{1} = \{u: \enskip (u,u)<\infty,\, (u',u')<\infty\}\,,
\end{equation*}
    and act as
\begin{equation*}
    Q_{a}(u) = - \kappa(a) |u(0)|^{2} + \int_{0}^{\infty} |u'(r)|^{2} dr 
\end{equation*}
    (here one can clearly see why the Friedrichs extension 
$ \kappa = -\infty $ for the set of forms
$ Q_{a}(u) $ bounded from below is called \emph{maximal}).
    The term ``extension'' also implies
    the extension condition
$ Q_{0} \subset Q_{a} $:
\begin{equation*}
    \Ww_{0}^{1} \subset \Ww_{(a)}^{1},\quad
	Q_{a}(u) = Q_{0}(u), \quad u\in\Ww_{0}^{1} .
\end{equation*}

    Thus we have obtained a one-to-one correspondence between the
    self-adjoint extensions
$ T_{a} $,
$ a\neq 0 $
    which act the same way as the second derivative operators
    but on different domains and the quadratic forms
$ Q_{a} $,
    defined on a single domain, but all acting differently.
    In the general case the structure of the domains of the quadratic forms
    corresponding to self-adjoined extensions of semi-bounded symmetric
    operator is described by the Krein theorem
\cite{Krein}.

    Applied to the case of the quadratic form of the Laplace operator in
    three dimensions
(\ref{SL}),
    the above equations imply that when closing the form starting
    from the domain of functions vanishing at zero, one arrives at
    a quadratic form on the domain of bounded functions.
    This form corresponds to the self-adjoint operator defined on the
    domain of bounded twice-differentiable functions (the ``common'' Laplace
    operator).
    But, at the same time, the original quadratic form
(\ref{SL}) can be extended to
    a closed form
$ Q_{a}^{\Delta}(f) $,
    which acts on functions diverging at zero as
$ r^{-1} $:
\begin{multline}
\label{Qfs}
    Q_{a}^{\Delta}(f) = \lim_{r\to 0} \Bigl(
= \sum_{j} \int_{\RR^{3}\setminus B_{r}}
    |\frac{\pl f}{\pl x_{j}}|^{2} d^{3} \vec{x}
- (\kappa(a)+\frac{1}{r}) \int_{\partial B_{r}} |f(\vec{x})|^{2} d^{2} s
\Bigr) =\\
=    \sum_{0\leq |m|\leq l} \int \bigl(|\frac{df_{lm}}{dr}|^{2}
	+ \frac{l(l+1)}{r^{2}}|f_{lm}|^{2}\bigr)r^{2}dr
- \kappa(a) \lim_{r\to0} \int_{\Sph^{2}}|f(r,\Omega)|^{2} r^{2} d\Omega ,
\end{multline}
    where
$ B_{r} $
    is a ball of radius
$ r $
    centered at the origin.
    The renormalization and scattering theory corresponding to
    such extensions were originally studied in
\cite{BF}.
    A thorough survey on perturbations and extension of differential
    operators can be found in
\cite{AK}.

\section{Vector Laplace operator}
\label{svector}
    Let us turn to the vector Laplace operator. This operator acts
    on 3-dimensi\-o\-nal vector functions
$ \vec{f}(\vec{x}) $ of 3 variables
    in the following way
\begin{equation*}
    \Delta \vec{f}(\vec{x}) = -\sum_{j=1}^{3} \frac{\pl^{2}}{\pl x_{j}^{2}}
	\vec{f}(\vec{x}) = \frac{1}{r^{2}} \Bigl( -\frac{\pl}{\pl r} r^{2}
	\frac{\pl}{\pl r} + \Dtph \Bigr) \vec{f}(\vec{x}(r,\Omega)) .
\end{equation*}
    It is natural to assume that in spherical
    coordinates the deficiency indices of a symmetric
    operator acting on the domain of functions vanishing at zero are
$ (3,3) $.

    Now, in place of the scalar spherical functions
$ Y_{lm} $, let us introduce
    three vector spherical harmonics (VSH)
\cite{VSH}:
\begin{align}
\label{VSH1}
    \vec{\Upsilon}_{lm} = & \frac{\vec{x}}{r} Y_{lm} , \quad
	0 \leq l, \quad |m| \leq l, \\
    \vec{\Psi}_{lm} = & (l(l+1))^{-1/2} r \vec{\pl} Y_{lm} , \quad
	1 \leq l , \quad |m| \leq l, \\
\label{VSH3}
    \vec{\Phi}_{lm} = & (l(l+1))^{-1/2} (\vec{x} \times \vec{\pl}) Y_{lm},
	\quad 1 \leq l , \quad |m| \leq l .
\end{align}
    Although VSH contain variable
$ r $
    in the definition, it is not hard to see that these functions
    are scale invariant and depend only on the angles
$ \Omega $.
    VSH defined as
(\ref{VSH1})--(\ref{VSH3})
    form an orthogonal and normalized system of vectors:
\begin{align}
    \int_{\Sph^{2}} \overline{\vec{\Upsilon}_{lm}(\Omega)}
        \vec{\Psi}_{l'm'}(\Omega) d\Omega & = 0 ,\quad
    \int_{\Sph^{2}} \overline{\vec{\Upsilon}_{lm}(\Omega)}
        \vec{\Upsilon}_{l'm'}(\Omega) d\Omega = \delta_{ll'} \delta_{mm'} , \\
    \int_{\Sph^{2}} \overline{\vec{\Upsilon}_{lm}(\Omega)}
        \vec{\Phi}_{l'm'}(\Omega) d\Omega       & = 0 ,\quad
    \int_{\Sph^{2}} \overline{\vec{\Psi}_{lm}(\Omega)}
        \vec{\Psi}_{l'm'}(\Omega) d\Omega = \delta_{ll'} \delta_{mm'} , \\
    \int_{\Sph^{2}} \overline{\vec{\Phi}_{lm}(\Omega)}
        \vec{\Psi}_{l'm'}(\Omega) d\Omega & = 0 ,\quad
    \int_{\Sph^{2}} \overline{\vec{\Phi}_{lm}(\Omega)}
        \vec{\Phi}_{l'm'}(\Omega) d\Omega = \delta_{ll'} \delta_{mm'} .
\end{align}
    Vector spherical harmonics allow for a unique representation
    of a vector function 
$ \vec{f}(\vec{x}) $
    as a series
\begin{equation}
\label{fext}
    \vec{f}(\vec{x}) =
	\sum_{0\leq |m| \leq l} y_{lm}(r) \vec{\Upsilon}_{lm} +
	\sum_{1 \leq l, |m| \leq l} \psi_{lm}(r) \vec{\Psi}_{lm} +
	\sum_{1 \leq l, |m| \leq l} \phi_{lm}(r) \vec{\Phi}_{lm} .
\end{equation}
    Although the separation of variables for the action of laplacian
$ \Delta $
    still can be written in the form
\begin{equation*}
    \Delta \bigl(z(r) \vec{Z}_{lm}\bigr) =
-\frac{1}{r^{2}} \frac{\pl}{\plr} r^{2} \frac{\pl}{\plr} z(r) \vec{Z}_{lm}
	+ \frac{z(r)}{r^{2}} \Dtph \vec{Z}_{lm}, \quad
	    \vec{Z} = \vec{\Upsilon}, \vec{\Psi}, \vec{\Phi} ,
\end{equation*}
    the harmonics are mixed in the second term.
    The action of
$ \Dtph $
    on VSH is not diagonal
(for $ l \geq 1 $),
    but with the normalization
(\ref{VSH1})--(\ref{VSH3})
    it appears to be symmetric:
\begin{align}
\nonumber
    \Dtph \vec{\Upsilon}_{lm} &= (2+l(l+1)) \vec{\Upsilon}_{lm}
	    - 2 \sqrt{l(l+1)} \vec{\Psi}_{lm} ,\\
\label{LVSH}
    \Dtph \vec{\Psi}_{lm} &= -2 \sqrt{l(l+1)} \vec{\Upsilon}_{lm}
	    + l(l+1) \vec{\Psi}_{lm} ,\\
\nonumber
    \Dtph \vec{\Phi}_{lm} &= l(l+1) \vec{\Phi}_{lm}
\end{align}
    (one may notice that the same expression holds 
    for the component
$ \vec{\Upsilon}_{00} $,
    for
$ l=0 $).
    Change of the basis
\begin{equation}
\label{Bcha}
    \begin{pmatrix} \vec{\Upsilon}_{lm} \\
	\vec{\Psi}_{lm}
    \end{pmatrix} \to \left(
    \begin{array}{l}
	\vec{\varUpsilon}_{lm} = (2l+1)^{-1/2} (\sqrt{l} \vec{\Upsilon}_{lm}
	    + \sqrt{l+1} \vec{\Psi}_{lm}) \\
	\vec{\varPsi}_{lm} = (2l+1)^{-1/2} (-\sqrt{l+1} \vec{\Upsilon}_{lm}
	    + \sqrt{l} \vec{\Psi}_{lm})
    \end{array} \right)
\end{equation}
    diagonalizes the action of
$ \Dtph $:
\begin{align*}
    \Dtph \vec{\varUpsilon}_{lm} & = (l-1)l \vec{\varUpsilon}_{lm} , \\
    \Dtph \vec{\varPsi}_{lm} & = (l+1)(l+2) \vec{\varPsi}_{lm} .
\end{align*}
    Here we may conclude that for
$ l=1 $
    the action of vector Laplace operator on the subspace of the
    component
$ \vec{\varUpsilon}_{1m} $
    coincides with the action of scalar Laplace operator
(\ref{TauL}) on the subspace
    of the spherical harmonics
\begin{equation*}
    \Delta \bigl(y(r) \vec{\varUpsilon}_{1m}\bigr) =
	- r^{-2}\frac{\pl}{\plr} r^{2} \frac{\pl}{\plr}y(r)
	    \vec{\varUpsilon}_{1m}, \quad m = -1, 0, 1. 
\end{equation*}
    Hence the vector Laplace operator
    written in basis
(\ref{VSH1})--(\ref{VSH3}),
    acting on smooth functions vanishing at zero
    also acquires the deficiency indices
$ (3,3) $.
    But now these numbers originate from the dimension of the subspace
    of spherical harmonics of momentum
$ l=1 $,
    rather than from the number of components of the vector field.

    The latter assertion was a small exercise. 
    From now on we forget the change
(\ref{Bcha})
    and turn to the transverse subspace.
    
\section{Transverse subspace}
\label{strav}
    Let us define the transverse subspace as a linear subset
    of vector functions satisfying
\begin{equation}
\label{transc}
    \vec{\pl} \cdot \vec{f}(\vec{x}) \equiv \sum_{j=1}^{3} \frac{\pl}{\pl x^{j}}
	f^{j}(\vec{x}) = 0 .
\end{equation}
    It is not hard to check that any product of
$ \vec{\Phi}_{lm} $
    by a radial function is transverse
\begin{multline*}
    \vec{\pl} \cdot \bigl(\phi(r) \vec{\Phi}_{lm} \bigr) = \vec{\pl} \phi(r)
	\cdot \vec{\Phi}_{lm} + \phi(r) \vec{\pl} \cdot \vec{\Phi}_{lm} = \\
    = (l(l+1))^{-1/2} \bigl(
    \phi'(r) r^{-1} \vec{x} \cdot (\vec{x}\times \vec{\pl}) Y_{lm}
    + \phi(r) \vec{\pl} \cdot (\vec{x}\times\vec{\pl}) Y_{lm}
	\bigr) = 0 ,
\end{multline*}
    here we have used the trivial properties of the mixed product
\begin{equation*}
    \vec{x} \cdot (\vec{x}\times\vec{\pl}) = 0 , \quad
    \vec{\pl} \cdot (\vec{x}\times\vec{\pl}) = 0 .
\end{equation*}
    The action of
$ \Delta $,
    along with that of its quadratic form, on the last sum of decomposition
(\ref{fext}),
    coincide with the action and the quadratic form of the scalar
    spherical laplacian
(\ref{Lsc}) for
$ l \geq 1 $.
    This operator is essentially self-adjoint on the subspaces corresponding to
    each $ l $.
    We are not going to describe it further, but rather turn to the second
    part of the transverse subspace.

    The first two sums of the decomposition
(\ref{fext})
    are not transverse when taken separately. 
    But with a special choice
    of the coefficient functions
\begin{align}
\label{upar1}
    y_{lm}(r) &= \sqrt{l(l+1)} \frac{u_{lm}(r)}{r^{2}} ,\\
\label{upar2}
    \psi_{lm}(r) &= \frac{u'_{lm}(r)}{r}
\end{align}
    each term 
$ y_{lm} \vec{\Upsilon}_{lm} + \psi_{lm} \vec{\Psi}_{lm} $
    becomes such:
\begin{align}
\label{treq}
    \vec{\pl} &\cdot \bigl(\sqrt{l(l+1)}\frac{u_{lm}}{r^{2}}\vec{\Upsilon}_{lm}
        +\frac{u'_{lm}}{r}\vec{\Psi}_{lm}\bigr) =\\
\nonumber
    &= \sqrt{l(l+1)} Y_{lm}
        \bigl( (\frac{u'_{lm}}{r^{2}}-\frac{2u_{lm}}{r^{3}})
        \frac{\vec{x}}{r}\cdot\frac{\vec{x}}{r} 
    + \frac{u_{lm}}{r^{2}} \vec{\pl}\cdot \frac{\vec{x}}{r} \bigr) +\\
\nonumber &\quad
    + 1/\sqrt{l(l+1)} u'_{lm} \vec{\pl}\cdot\vec{\pl} Y_{lm} = 0 .
\end{align}
    Here we have taken into account the relations
\begin{equation*}
    \vec{x}\cdot\vec{\pl}Y_{lm}(\Omega) = 0 , \quad
    \vec{\pl} \cdot \frac{\vec{x}}{r} = \frac{2}{r} , \quad
    \vec{\pl}\cdot\vec{\pl} Y_{lm} = -\frac{l(l+1)}{r^{2}} Y_{lm} .
\end{equation*}
    Equality
(\ref{treq})
    also immediately follows from the representation
\begin{equation*}
    \sqrt{l(l+1)}\frac{u}{r^{2}}\vec{\Upsilon}_{lm}
        + \frac{u'}{r}\vec{\Psi}_{lm} = \vec{\partial}\times\frac{u}{r}
    \vec{\Phi}_{lm} ,
\end{equation*}
    prove of which, however, requires some two-line calculation.
    The choice of parameter
$ u_{lm} $
    in
(\ref{upar1})--(\ref{upar2}) 
    is definitely not unique, and we use this particular one
    to obtain operators similar to
(\ref{Tl}).

    Let us take the functions
$ u_{lm} $
    to be smooth
    on the positive half-axis and vanishing at zero together
    with their derivatives. We denote this subset as
$ \SP_{l}^{0} $
    and the subset of the transverse functions generated by such
$ u_{lm} $
    as
$ \DD_{lm} $:
\begin{equation*}
    \DD_{lm} = \bigl\{ \vec{f} =
	\sqrt{l(l+1)}\frac{u_{lm}(r)}{r^{2}}\vec{\Upsilon}_{lm}(\Omega)
        +\frac{u'_{lm}(r)}{r}\vec{\Psi}_{lm}(\Omega),
	    \quad u_{lm} \in\SP_{l}^{0}
    \bigr\} .
\end{equation*}
    The calculation of the action of the Laplace operator on
$ \DD_{lm} $ by virtue of
(\ref{LVSH})
\begin{multline*}
    \Delta \bigl(\sqrt{l(l+1)}\frac{u_{lm}}{r^{2}}\vec{\Upsilon}_{lm}
	+\frac{u'_{lm}}{r}\vec{\Psi}_{lm}\bigr) =\\
    = \sqrt{l(l+1)}\frac{1}{r^{2}} \bigl(-u''_{lm}
	+\frac{l(l+1)}{r^{2}}u_{lm}\bigr) \vec{\Upsilon}_{lm} 
    + \frac{1}{r} \bigl(-u''_{lm} +\frac{l(l+1)}{r^{2}}u_{lm}\bigr)'
	\vec{\Psi}_{lm} 
\end{multline*}
    reveals that, first -- the subset
$ \DD_{lm} $
    is an invariant subspace for
$ \Delta $,
    and second -- that the laplacian induces on functions
$ u_{lm} $
    from
$ \SP_{l}^{0} $
    the action of second order differential operator
$ T_{l} $:
\begin{equation}
\label{Lact}
    \Delta : \quad u_{lm} \to T_{l} u_{lm} 
	= - \frac{\pl^{2}}{\plr^{2}}u_{lm} + \frac{l(l+1)}{r^{2}} u_{lm} .
\end{equation}

    The scalar product
(\ref{fgprod})
    from the space of 3-dimensional vector functions,
    being contracted to
$ \DD_{lm} $,
    also naturally transfers to
    sesquilinear form on
$ \SP_{l}^{0} $:
\begin{multline}
    \int \overline{\vec{f}_{v}(\vec{x})} \vec{f}_{u}(\vec{x}) d^{3}x =\\
    =\int
    \bigl(\sqrt{l(l+1)} \frac{\ol{v_{lm}}}{r^{2}}\ol{\vec{\Upsilon}_{lm}}
	+\frac{\ol{v'_{lm}}}{r}\ol{\vec{\Psi}_{lm}} \bigr) \cdot
    \bigl(\sqrt{l(l+1)} \frac{u_{lm}}{r^{2}}\vec{\Upsilon}_{lm}
	+\frac{u'_{lm}}{r} \vec{\Psi}_{lm} \bigr) r^{2} dr d\Omega =\\
\label{ff}
    = \int \bigl( \ol{v'_{lm}}u'_{lm} + \frac{l(l+1)}{r^{2}}\ol{v_{lm}}u_{lm}
    \bigr) dr \equiv \langle v,u \rangle_{l} .
\end{multline}
    In this way each linear subset
$ \SP_{l}^{0} $
    acquires a Hilbert structure.
    Along with the action of
$ T_{l} $
    this structure defines a quadratic form in the space
$ \SP_{l}^{0} $
    and the action of this form corresponds to the action of the quadratic
    form of vector laplacian on the subspace $ \DD_{lm} $.

\section{Operator of the quadratic form}
\label{soper}
    We are not going to construct a closed operator out of
$ T_{l} $
    in the sesquilinear form 
(\ref{fgprod})
    as well as the closed quadratic form
\begin{equation}
\label{forml}
    (\vec{f}_{u},\Delta \vec{f}_{u})
        = \langle u, T_{l}u\rangle_{l}, \quad l\geq 1
\end{equation}
    in that scalar product.
    Instead we shall consider a somewhat more general problem.
    Namely, the extensions of the quadratic form
(\ref{forml})
    with respect to the scalar product
(\ref{plainprod}).
    The reason for this is the consideration that
    underlying physics, described by the quadratic form of the
    transverse laplacian, does not necessarily restrict its states by
    normalization
(\ref{ff}).
    Another reason is that 
(\ref{plainprod}) is
    one of products that are local in variable
$ r $
    (as we mentioned, one could freely transfer part of the operator into
    the product, but usually this would result in a nonlocal expression).

    By means of integration by parts one can see that for
$ u, v \in \SP_{l}^{0} $
    the following equality holds
\begin{align}
\nonumber
    \langle v,&T_{l} u \rangle_{l} = \\
\nonumber
&= \int_{0}^{\infty} \bigl(\frac{d \bar{v}}{dr} \frac{d}{dr}
(-\frac{d^{2} u}{dr^{2}}+\frac{l(l+1)}{r^{2}}u) + \frac{l(l+1)}{r^{2}}\bar{v}
(-\frac{d^{2} u}{dr^{2}}+\frac{l(l+1)}{r^{2}}u) \bigr) dr =\\
\nonumber
&= \int_{0}^{\infty}
(-\frac{d^{2} \bar{v}}{dr^{2}}+\frac{l(l+1)}{r^{2}}\bar{v})
(-\frac{d^{2} u}{dr^{2}}+\frac{l(l+1)}{r^{2}}u) dr = (T_{l}v, T_{l}u) \\
\nonumber
& = \int_{0}^{\infty} \bar{v} \bigl(-\frac{d^{2}}{dr^{2}} 
(-\frac{d^{2} u}{dr^{2}}+\frac{l(l+1)}{r^{2}}u) + \frac{l(l+1)}{r^{2}}
(-\frac{d^{2} u}{dr^{2}}+\frac{l(l+1)}{r^{2}}u) \bigr) dr =\\
\label{vTu}
    & = (v, T_{l}^{2} u) .
\end{align}
    This implies that quadratic form
$ \langle u , T_{l} u \rangle_{l} $
    in the scalar product
(\ref{plainprod})
    is defined by the fourth order differential operator
\begin{equation}
\label{Tl2}
    T_{l}^{2} = \bigl(-\frac{d^{2}}{dr^{2}} + \frac{l(l+1)}{r^{2}}\bigr)^{2} .
\end{equation}

\subsection{Alternative interpretation of the quadratic form}
    We can also provide an alternative interpretation of the possible
    extensions of the operator
$ T_{l} $,
$ l=1 $
    in the product
(\ref{ff}). Let
$ E_{l} $
    be a linear differential operation acting on functions from
$ \SP_{0}^{l} $
    in the following way
\begin{equation*}
    E_{l} u = r^{-l} \frac{d}{dr} (r^{l}u) .
\end{equation*}
    Adjoint (in the sence of integration by parts in the scalar product)
    operation
$ E_{l}^{*} $
    looks as
\begin{equation*}
    E_{l}^{*} u = - r^{l} \frac{d}{dr} (r^{-l}u) = E_{-l} u,
\end{equation*}
    and the products of
$ E_{l} $ and
$ E_{l}^{*} $
    are equal to the described second order differential operations
$ T_{l} $ and
$ T_{l-1} $:
\begin{align*}
    E_{l}^{*} E_{l}u = & -r^{l} \frac{d}{dr} r^{-2l} \frac{d}{dr} (r^{l}u)
        = T_{l} \\
    E_{l} E_{l}^{*}u = & -r^{-l} \frac{d}{dr} r^{2l} \frac{d}{dr} (r^{-l}u)
        = T_{-l}u = T_{l-1}u .
\end{align*}
    For any
$ u,v \in \SP_{l}^{0} $
    the scalar product
(\ref{ff})
    takes the form
\begin{equation*}
    \langle u, v \rangle_{l} = (E_{l}u, E_{l}v) ,
\end{equation*}
    while the form
(\ref{vTu})
    is represented by the expression
\begin{equation*}
    \langle u, T_{l} v \rangle_{l} = (E_{l}u, E_{l} T_{l}v)
        = (E_{l}u, E_{l} E_{l}^{*} E_{l}v) .
\end{equation*}
    The linear map
\begin{equation*}
    u \to \psi_{u} = E_{l} u = r^{-l} \frac{d}{dr} r^{l} u(r)
\end{equation*}
    has a trivial kernel on
$ \SP_{l}^{0} $,
    and naturally transfers the sesquilinear form
(\ref{ff})
    into
(\ref{plainprod})
\begin{equation*}
    \langle u, v \rangle_{l} = (\psi_{u},\psi_{v}) .
\end{equation*}
    In the terms of
$ \psi $
    the quadratic form
(\ref{forml})
    looks as simple as
\begin{equation}
\label{formpsi}
    \langle u, T_{l} v \rangle_{l} = (E_{l}u, E_{l} E_{l}^{*} E_{l}v)
        = (\psi_{u}, E_{l} E_{l}^{*} \psi_{u}) = (\psi_{u}, T_{l-1}\psi_{u}) .
\end{equation}
    Thus the problem transfers to the construction of the extensions of
    quadratic forms of operators
$ T_{l-1} $
    in the ``plain'' scalar product
(\ref{plainprod}).
    These operators, as was said in
\ref{saradial},
    are essentially self-adjoint for
$ l>1 $
    and are symmetric with deficiency indices
$ (1,1) $
    for
$ l=1 $.

\section{Deficiency indices of $ T_{l}^{2} $}
\label{sdef}
    The solutions of the equation
\begin{equation*}
    - \frac{d^{2}g}{dr^{2}} + \frac{l(l+1)}{r^{2}} g = \lambda g
\end{equation*}
    coincide, up to the factor of
$ r $,
    with the spherical Bessel functions. Hence, in order to describe
    the solutions of an equation with
$ T_{l}^{2} $
    in the LHS
    we need to provide some basic facts from this classical theory.

\subsection{Some properties of the spherical Bessel functions}
    Let
$ D_{l} $ be a linear differential operation defined as
\begin{equation}
\label{Dw}
    D_{l} w(r) = r^{l+1} \bigl(\frac{1}{r}\frac{d}{dr}\bigr)^{l} \frac{w}{r}
	= r^{l}\bigl(\frac{d}{dr}\frac{1}{r}\bigr)^{l} w(r) .
\end{equation}
    Then, via mathematical induction and some differentiation one
    obtains the equality
\begin{equation}
\label{QD}
    T_{l} D_{l} w(r) = D_{l} T_{0} w(r) = -D_{l} \frac{d^{2}w}{dr^{2}} .
\end{equation}
    This allows one to construct eigenfunctions of the operator
$ T_{l}^{2} $
    from the eigenfunctions of
\begin{equation*}
    T_{0}^{2} = \frac{d^{4}}{dx^{4}} .
\end{equation*}

    The behaviour of the functions represented as
(\ref{Dw})
    in the vicinity of zero can be described as follows. Let
\begin{equation*}
    w(r) = r^{k} ,
\end{equation*}   
    then
\begin{equation*}
    D_{l} w(r) = (k-2l+1) \ldots (k-1) r^{k-l} .
\end{equation*}
    In particular, if
$ w(r) $ has an expansion near zero which looks like
\begin{equation*}
    w(r) = w_{0} + w_{1} r + \ldots + w_{5} r^{5}
	+\OO(r^{6}),
\end{equation*}
    then
\begin{align}
\label{A1}
    D_{1} w(r) &= -\frac{w_{0}}{r} + w_{2}r + 2w_{3}r^{2} + 3w_{4} r^{3}
	+ \OO(r^{4}),\\
\label{A2}
    D_{2} w(r) &= \frac{3w_{0}}{r^{2}} -w_{2} + 3w_{4}r^{2} +8w_{5} r^{3}
	+ \OO(r^{4}) .
\end{align}
    That is, the operation
$ D_{l} $
    lowers the power of its argument by 
$ l $,
    and simultaneously
$ l $
    powers are eliminated.
    This also means that one has to separately check that after applying
$ D_{l} $
    to a set of independent functions the resulting set is still independent.
    We will be dealing with exponents of non-commuting periods
    for which this fact is rather trivial, so we will not focus on this. 

\subsection{Frobenius method}
\label{Frob}
    The deficiency indices of symmetric operator
$ T_{l} $,
    that is, the dimensions of the subspaces of square integrable solutions
    of the equations
\begin{equation}
\label{Q}
    \bigl(-\frac{d^{2}}{dr^{2}}+\frac{l(l+1)}{r^{2}}\bigr) g
	= \pm i\tilde{\rho}^{2} g
\end{equation}
    can be estimated by the Frobenius method
    (see {\it e.g.}
\cite{Richt})
    of expansion of
$ g(r) $
    in a vicinity of a critical point (in our case it is point zero).
    In the case of the second order differential operator
$ T_{l} $
    this method results in the following equation for the power
$ \alpha $
\begin{equation}
\label{alpha}
    -\frac{d^{2} r^{\alpha}}{dr^{2}} + \frac{l(l+1)}{r^{2}} r^{\alpha} = 0 .
\end{equation}
    It gives two possible solutions of
(\ref{Q})
    behaving as 
\begin{equation*}
    g_{1}(r) \stackrel{r\to 0}{\sim} r^{l+1} , \quad
    g_{2}(r) \stackrel{r\to 0}{\sim} r^{-l} .
\end{equation*}
    The first solution is regular at zero but for
$ \RE\, \tilde{\rho}^{2} \,>\, 0 $
    it is exponentially growing at infinity.
    The second, in contrast, is regular at infinity, but
    for
$ l \,>\, 0 $
    it diverges at zero.

    A simple modification of the Frobenius method reveals, that
    two of the four possible solutions of the equation
\begin{equation}
\label{QQ}
\bigl(-\frac{d^{2}}{dr^{2}}+\frac{l(l+1)}{r^{2}}\bigr)^{2} g = \pm i\rho^{4} g
\end{equation}
    have the same behaviour at zero as the solutions of
(\ref{Q})
    (more precisely, they coincide with those of
(\ref{Q}) after the substitution
$ i\rho^{4} = \tilde{\rho}^{4} $),
    while the other two have powers of $ r $ near zero greater by two:
\begin{equation*}
    g_{3}(r) \stackrel{r\to 0}{\sim} r^{l+3} , \quad
    g_{4}(r) \stackrel{r\to 0}{\sim} r^{-l+2} 
\end{equation*}
    (in this case the first operator 
$ T_{l} $ in
(\ref{QQ})
    lowers the power by two, while the second one solves
(\ref{alpha})).
    The solutions
$ g_{1} $ and
$ g_{3} $
    grow exponentially at infinity,
$ g_{2} $
    is still not square-integrable at zero for
$ l \geq 1 $,
    while
$ g_{4} $
    does satisfy the square integrability if
$ l=1,2 $.
    That is, one can expect that the symmetric semi-bounded operator    
$ T_{l}^{2} $
    has the deficiency indices
$ (1,1) $
    when
$ l=1,2 $.

\subsection{Solutions of the fourth order differential equations}
    Let us construct the solutions of
(\ref{QQ})
    in an explicit form. Equation
(\ref{QD}),
    applied to the exponents
\begin{equation*}
    w = \exp\{e^{\pm i\frac{\pi}{8} \pm i\frac{\pi k}{2}} \rho r\} ,
	\quad k=0,1,2,3\,,
\end{equation*}
    gives
\begin{multline*}
    T_{l}^{2} D_{l} 
    \exp\{e^{\pm i\frac{\pi}{8} \pm i\frac{\pi k}{2}} \rho r\} =
    D_{l} T_{0}^{2}
    \exp\{e^{\pm i\frac{\pi}{8} \pm i\frac{\pi k}{2}} \rho r\} = \\
=    \pm i \rho^{4} D_{l} 
    \exp\{e^{\pm i\frac{\pi}{8} \pm i\frac{\pi k}{2}} \rho r\} .
\end{multline*}
    The functions
$ D_{l} \exp\{e^{\pm i\frac{\pi}{8} \pm i\frac{\pi k}{2}} \rho r\} $
    exponentially grow at infinity if
$ k = 0,3 $
    and exponentially vanish if
$ k = 1,2 $.
    The asymptotics
(\ref{A1})
    and
(\ref{A2})
    show, that within the linear span of exponents
$ k=1,2 $
    only the combinations proportional to their difference are
    square-integrable at zero (it is only in this difference
    where the coefficient
$ w_{0} $ vanishes).

    Thus we can see, that the subspaces of square integrable solutions of
(\ref{QQ})
    are one-dimensional and generated by the functions
\begin{align*}
    g_{l+} &= D_{l}\bigl(
	\exp\{e^{i\frac{5\pi}{8}} \rho r\}-\exp\{e^{i\frac{9\pi}{8}} \rho r\} 
    \bigr) , \\
    g_{l-} &= D_{l}\bigl(
	\exp\{e^{i\frac{7\pi}{8}} \rho r\}-\exp\{e^{i\frac{11\pi}{8}} \rho r\} 
    \bigr) .
\end{align*}

\section{Extensions of symmetric operator
$ T_{l}^{2} $}
\label{ssymm}
    The functions
$ g_{l\pm} $
    constructed in the previous section generate kernels of the corresponding
    adjoint operators:
\begin{equation*}
    (T_{l}^{2} + i\rho^{4})^{*} g_{l+} = 0 , \quad
    (T_{l}^{2} - i\rho^{4})^{*} g_{l-} = 0 .
\end{equation*}
    These kernels are the only orthogonal complements of the ranges
    of the shifted symmetric operators
\begin{equation*}
    Ran_{\pm} = \{ g = (T_{l}^{2} \pm i\rho^{4}) u , \quad u \in \SP_{l}^{0}\} .
\end{equation*}
    The Cayley transform provides every symmetric operator
$ T_{l}^{2} $
    with a partial isometry
$ U $,
    mapping from
$ Ran_{+} $ to
$ Ran_{-} $
    according to the rule
\begin{equation*}
    U : (T_{l}^{2} +i\rho^{4}) u \to (T_{l}^{2} -i\rho^{4}) u .
\end{equation*}
    Based on the above description of the kernels, the isometric operator
$ U $
    can be extended to an unitary operator
$ U_{a} $
    by letting the latter act on the orthogonal
    complements of
$ Ran_{+} $
    as
\begin{equation*}
    U_{a} e^{ia} g_{l+} = e^{-ia} g_{l-} , \quad 0\leq a <\pi ,
\end{equation*}
    where
$ e^{2ia} $ is an unitary parameter.
    The operator
$ U_{a} $
    is a Cayley transform of some self-adjoint extension
$ T_{la}^{2} $
    of the symmetric operator
$ T_{l}^{2} $.
    The domain of 
$ T_{la}^{2} $
    now contains the vector
$ h_{l} $
    which satisfies the conditions
\begin{align*}
    (T_{la}^{2} + i\rho^{4}) h_{l}^{a} &= e^{ia} g_{l+} , \\
    (T_{la}^{2} - i\rho^{4}) h_{l}^{a} &= e^{-ia} g_{l-} .
\end{align*}
    Choosing
$ h_{l}^{a} $
    to be a linear combination of
$ g_{l+} $ and
$ g_{l-} $
    one can find 
\begin{multline*}
    h_{l}^{a}
    = \frac{e^{ia}}{2i\rho^{4}} g_{l+} - \frac{e^{-ia}}{2i\rho^{4}}g_{l-}
    = \frac{1}{2i\rho^{4}} D_{l} \bigl(
\exp\{ia + e^{i\frac{5\pi}{8}}\rho r\} - \exp\{ia - e^{i\frac{\pi}{8}}\rho r\}
    -\\
-\exp\{-ia - e^{-i\frac{\pi}{8}}\rho r\}+\exp\{-ia +e^{-i\frac{5\pi}{8}}\rho r\}
    \bigr).
\end{multline*}
    Thus we arrive at a conclusion that the essential domain of
$ T_{la}^{2} $
    is a sum
\begin{equation*}
    \SP_{l}^{a} = \SP_{l}^{0} \dotplus h_{l}^{a},
\end{equation*}
    and that the action of
$ T_{la}^{2} $
    is the fourth order differential operation
(\ref{Tl2}).

    Using the asymptotics
(\ref{A1}), (\ref{A2})
    one can calculate the first coefficients of the expansion of
$ h_{l}^{a} $
    in the vicinity of zero
\begin{align*}
    h_{1}^{a}(r) & = \sin a \, \rho r
	+ \frac{\sqrt{2}}{3}\cos(a+\frac{\pi}{8})\, \rho^{3} r^{2}
	+ \OO(r^{4}) \\
    h_{2}^{a}(r) & = \sin a \, \rho^{2}
	+ \frac{\sqrt{2}}{15}\cos(a-\frac{\pi}{8})\, \rho^{5} r^{3}
	+ \OO(r^{4}) .
\end{align*}
    These expansions yield the boundary conditions for the domains of
$ T_{la}^{2} $:
\begin{align}
\label{B1}
    l=1: \quad h_{1}^{a}{}''(0) &= \frac{2\sqrt{2}}{3}
    \frac{\cos(a+\pi/8)}{\sin a}
	\rho h_{1}^{a}{}'(0) , \quad h_{1}^{a}(0) = h_{1}^{a}{}'''(0) = 0 , \\
\label{B2}
    l=2: \quad h_{2}^{a}{}'''(0) &= \frac{2\sqrt{2}}{5}
	\frac{\cos(a-\pi/8)}{\sin a}
	\rho^{3} h_{1}^{a}(0) , \quad h_{1}^{a}{}'(0) = h_{1}^{a}{}''(0) = 0 .
\end{align}
    Hence the essential domain
$ \SP_{l}^{a} $
    can be described as a linear set of functions bounded by the norm
\begin{equation*}
    \bigl(\: \cdot \:, \:\cdot\: \bigr) +
    \bigl( (-\frac{d^{2}}{dr^{2}} +\frac{l(l+1)}{r^{2}})^{2}\:\cdot \: ,
	(-\frac{d^{2}}{dr^{2}} +\frac{l(l+1)}{r^{2}})^{2}\: \cdot\: \bigr)
\end{equation*}
    and satisfying the boundary conditions
(\ref{B1}),
(\ref{B2}), correspondingly.

\subsection{Integration by parts and the property of symmetry}
    Now let us check that the conditions
(\ref{B1}) and (\ref{B2})
    are in agreement with the symmetry property of the operators
$ T_{la}^{2} $:
\begin{multline*}
    \int \bar{v}
\bigl(-\frac{d^{2}}{dr^{2}} +\frac{l(l+1)}{r^{2}}\bigr)^{2} u dr 
    -\int \bigl(-\frac{d^{2}}{dr^{2}}
	+\frac{l(l+1)}{r^{2}}\bigr)^{2}\bar{v} u dr = \\
= \bigl(\bar{v}u''' -\bar{v}'''u -\bar{v}'u'' +\bar{v}''u'
	+2\frac{l(l+1)}{r^{2}}(\bar{v}'u-\bar{v}u')\bigr)\bigr|_{0}^{\infty}.
\end{multline*}
    The upper limit in the RHS is zero for the reason that
$ u $,
$ v $,
    as well as their derivatives, vanish at infinity.
    The lower limit can be calculated by expanding the expression
$ (\bar{v}'u-\bar{v}u') $
    in the values of
$ u $ and
$ v $
    and their derivatives at zero (we denote them as
$ u_{0} $, $ u_{0}' $, $ u_{0}'' $, $ u_{0}''' $ and similarly for 
$ v $).
    Then RHS above transforms into three terms
\begin{multline*}
    2\frac{l(l+1)}{r^{2}}(\bar{v}_{0}u_{0}'-\bar{v}_{0}'u_{0})\bigr|_{r=0}
    +2\frac{l(l+1)}{r}(\bar{v}_{0}u_{0}''-\bar{v}_{0}''u_{0})\bigr|_{r=0} +\\
    +(l(l+1)-1)(\bar{v}_{0}u_{0}''' -\bar{v}_{0}'''u_{0}) 
    +(l(l+1)+1)(\bar{v}_{0}'u_{0}'' -\bar{v}_{0}''u_{0}') .
\end{multline*}
    In order for this sum to be equal to zero each of these terms
    must be zeroed separately. This leads to the conditions
\begin{equation}
\label{CC}
    \frac{u_{0}'}{u_{0}} = \frac{\bar{v}_{0}'}{\bar{v}_{0}} , \quad
    \frac{u_{0}''}{u_{0}} = \frac{\bar{v}_{0}''}{\bar{v}_{0}} , \quad
    \frac{u_{0}'''}{u_{0}} = \frac{\bar{v}_{0}'''}{\bar{v}_{0}} , 
\end{equation}
    which uniquely bind three of the first four coefficients of the expansion of
$ u $ (and $ v $)
    to the fourth one.

    It is not hard to see, that the conditions
(\ref{B1})
    and
(\ref{B2})
    are in agreement with
(\ref{CC}),
    and in this way the differential operators
$ T_{la}^{2} $
    are symmetric on the domains
$ \SP_{l}^{a} $.

\section{Discrete spectrum}
\label{sdiscrete}
    Using the boundary conditions
(\ref{B1})
    and
(\ref{B2})
    let us try to study, for which values of
$ a $
    there exist (discrete) eigenvalues of 
$ T_{la}^{2} $.
    We shall restrict ourselves to the negative spectrum, while the absence
    of positive eigenvalues will be discussed in the next section.
    Suppose that function
$ \tv_{l}^{\kappa} $
    satisfies the equation
\begin{equation*}
    \bigl(-\frac{d^{2}}{dr^{2}}+\frac{l(l+1)}{r^{2}}\bigr)^{2}
\tv_{l}^{\kappa}(r) = -\kappa^{4} \tv_{l}^{\kappa}(r) , \quad \kappa >0.
\end{equation*}
    Let us write this function as
\begin{equation*}
    \tv_{l}^{\kappa}(r) = D_{l} w^{\kappa}(r) ,
\end{equation*}
    then
\begin{equation*}
    T_{l}^{2} \tv_{l}^{\kappa} = T_{l}^{2} D_{l} w^{\kappa} 
	= D_{l} T_{0}^{2} w^{\kappa} = D_{l} \frac{d^{4}}{dr^{4}} w^{\kappa}
	= -\kappa^{4} D_{l} w^{\kappa} .
\end{equation*}
    This equation shows that
$ w^{\kappa} $,
    up to a solution of
\begin{equation*}
    D_{l}w = 0 ,
\end{equation*}
    belongs to the linear span of four exponents
\begin{equation*}
    \exp\{e^{i\frac{\pi}{4}+i\frac{\pi k}{2}} \kappa r \} ,
	\quad k = 0, 1, 2, 3 .
\end{equation*}
    Two of these exponents, namely those with
$ k=0,3 $,
    grow at infinity, while among the other two
    a square integrable solution for
$ \tv_{l}^{\kappa} $
    is produced only by the difference
\begin{equation*}
    w^{\kappa} = i\exp\{e^{-i\frac{3\pi}{4}} \kappa r \} -
    i\exp\{e^{i\frac{3\pi}{4}} \kappa r \}
\end{equation*}
    (the imaginary coefficients here are chosen in order to obtain a real
    value for 
$ \tv_{l}^{\kappa} $).
    Using the expansion of
$ w^{\kappa} $
    in the vicinity of zero
\begin{equation*}
    w^{\kappa} = \sqrt{2}\kappa r - \kappa^{2} r^{2}
	+ \frac{\sqrt{2}}{6}\kappa^{3} r^{3}
	- \frac{\sqrt{2}}{120}\kappa^{5} r^{5} + \OO(r^{6}) ,
\end{equation*}
    and the asymtotics
(\ref{A1})
    and
(\ref{A2})
    we find the first coefficients of
$ \tv_{l}^{\kappa} $:
\begin{align}
\label{V1}
    \tv_{1}^{\kappa} &= -\kappa^{2} r +\frac{\sqrt{2}}{3} \kappa^{3} r^{2}
	+ \OO(r^{4}) ,\\
\label{V2}
    \tv_{2}^{\kappa} &= \kappa^{2} -\frac{\sqrt{2}}{15} \kappa^{5} r^{3}
	+ \OO(r^{4}) .
\end{align}
    In order for the functions
$ v_{l}^{\kappa} $
    to be in the domains of
$ T_{la}^{2} $
    they should satisfy the corresponding boundary conditions
(\ref{B1}) or
(\ref{B2}).
    The second parts of these conditions are obviously obeyed because of the absence
    of the corresponding coefficients in
(\ref{V1})
    or
(\ref{V2}).
    Meanwhile, the first parts produce conditions which express
$ \kappa $
    in terms of the parameter
$ a $:
\begin{align}
\label{ka1}
    l=1: \quad & \kappa = -\rho \frac{\cos(a+\pi/8)}{\sin a} ,
	\quad 0\leq a < \pi , \\
\label{ka2}
    l=2: \quad & \kappa^{3} = -\rho^{3} \frac{\cos(a-\pi/8)}{\sin a} ,
	\quad 0\leq a < \pi .
\end{align}
    Each of the RHS of these equations is an one-to-one correspondence relating the interval
$ 0 \leq a < \pi $
    to the real axis joined with point
$ -\infty $.
    That is, there is a one-to-one correspondence between 
$ a $ and
$ \kappa $
    (at fixed $ \rho $).
    In order to simplify further calculations we will use 
    the dimensional parameter
$ \kappa $
    instead of both
$ a $ and
$ \rho $ without regard for the existence of discrete spectrum.
    With this in mind, the boundary conditions
(\ref{B1})
    and
(\ref{B2})
    take the form
\begin{align}
\label{NB1}
    l=1: \quad v_{1}^{\kappa}{}''(0) &= -\frac{2\sqrt{2}}{3} \kappa
	v_{1}^{\kappa}{}'(0) , \quad v_{1}^{\kappa}(0)
	= v_{1}^{\kappa}{}'''(0) = 0 , \\
\label{NB2}
    l=2: \quad v_{2}^{\kappa}{}'''(0) &= -\frac{2\sqrt{2}}{5} \kappa^{3}
	v_{2}^{\kappa}(0) , \quad v_{2}^{\kappa}{}'(0)
	= v_{2}^{\kappa}{}''(0) = 0 .
\end{align}
    Equations
(\ref{ka1})
    and
(\ref{ka2})
    show, that positive values of
$ \kappa $
    (which are equivalent to the existence of the negative discrete spectrum)
    correspond to parameter
$ a $
    residing in the intervals
\begin{align*}
    l = 1 : \quad & \frac{3\pi}{8} < a < \pi ,\\
    l = 2 : \quad & \frac{5\pi}{8} < a < \pi .
\end{align*}
    Thus we may conclude that 
    in the above intervals of
$ a $
    the self-adjoint operators
$ T_{la}^{2} $
    acquire negative eigenvalues
$ -\kappa^{4} $
    of multiplicity one,
    corresponding to (unnormalized) eigenfunctions
\begin{equation}
\label{tv}
    \tv_{l}^{\kappa} = iD_{l} \bigl(
        \exp\{e^{-i\frac{3\pi}{4}} \kappa r \} -
	\exp\{e^{i\frac{3\pi}{4}} \kappa r \} \bigr) .
\end{equation}

\section{Continuous spectrum}
\label{scont}
    The description of the continuous spectrum of operator
$ T_{l\kappa}^{2} $,
$ \kappa=\kappa(a) $
    involves the solution of the equation
\begin{equation*}
    \bigl(-\frac{d^{2}}{dr^{2}} +\frac{l(l+1)}{r^{2}}\bigr)^{2}
	\tu_{l}^{\lambda}
	= \lambda^{4} \tu_{l}^{\lambda} , \quad \lambda > 0 .
\end{equation*}
    As in the previous section, one can substitute
\begin{equation}
\label{uDw}
    \tu_{l}^{\lambda} = D_{l} w^{\lambda} ,
\end{equation}
    and then a simple fourth order equation arises for
$ w^{\lambda} $
\begin{equation*}
    \frac{d^{4}w^{\lambda}}{dr^{4}} = \lambda^{4} w^{\lambda} .
\end{equation*}
    Its solutions are the four linearly independent exponents
\begin{equation*}
    \exp\{i^{k} \lambda r\} , \quad k=0,1,2,3,
\end{equation*}
    which allow us to construct two linearly independent functions vanishing
    at zero and nongrowing at infinity.
    Any linear combination of these functions can be written
    (up to a common factor) as
\begin{equation*}
    w^{\lambda} = \sin\lambda r
	+ \sigma(\lambda) (\cos\lambda r -e^{-\lambda r}) .
\end{equation*}
    This expression immediately yields the absence of the positive
    discrete spectrum of
$ T_{la}^{2} $, since no any
$ w^{\lambda} $
    vanishes at infinity, nor does the main term in
$ \tu_{l}^{\lambda} $.

    Let us show that for any
$ \lambda > 0 $
    there is a single
$ \sigma(\lambda) $
    such that the function
(\ref{uDw})
    obeys the boundary condition
(\ref{NB1})
    or
(\ref{NB2}).
    Indeed, the expansion of
$ w^{\lambda} $
    in
$ r $
    near zero reads as
\begin{equation*}
    w^{\lambda}(r) = (1+\sigma) \lambda r - \sigma \lambda^{2} r^{2}
	+\frac{\sigma -1}{6}\lambda^{3} r^{3}
	+\frac{\sigma +1}{120} \lambda^{5} r^{5} + \OO(r^{6}) .
\end{equation*}
    Then the asymptotics 
(\ref{A1})
    and
(\ref{A2})
    give the estimates
\begin{align*}
    \tu_{1}^{\lambda} &= D_{1} w^{\lambda} = -\sigma\lambda^{2} r
	+ \frac{\sigma -1}{3} \lambda^{3} r^{2} + \OO(r^{4}) , \\
    \tu_{2}^{\lambda} &= D_{2} w^{\lambda} = \sigma\lambda^{2} 
	+ \frac{\sigma +1}{15} \lambda^{5} r^{3} + \OO(r^{4}) .
\end{align*}
    By comparing these expansions with
(\ref{NB1})
    and
(\ref{NB2})
    one may conclude that the dependence of
$ \sigma $
    on 
$ \lambda $
    is single valued and looks as 
\begin{align}
\label{S1}
    l=1: \quad
    \sigma_{1}(\lambda) &= \frac{\lambda}{\lambda -\sqrt{2}\kappa} ,
	\quad \lambda >0 ,\\
\label{S2}
    l=2: \quad
    \sigma_{2}(\lambda) &= \frac{\lambda^{3}}{\sqrt{2}\kappa^{3}-\lambda^{3}} ,
	\quad \lambda >0 .
\end{align}
    Thus we arrive at a continuous spectrum of multiplicity one,
    occupying all values on the positive half-axis.

\section{Normalization of eigenfunctions}
\label{snorm}
    To calculate the normalization of the eigenfunctions
$ \tv_{l}^{\kappa} $
    and of the ``eigenfunctions of the continuous spectrum''
$ \tu_{l}^{\lambda} $
    it is convenient to exploit the following identities for
    the operations
$ D_{l} $
\begin{align}
\label{N1}
    \int_{0}^{\infty} D_{1} \tw(r) D_{1} w(r) dr &=
	\int_{0}^{\infty} \frac{d\tw}{dr} \frac{dw}{dr} dr
	    - \frac{\tw w}{r} \Bigr|_{0}^{\infty} ,\\
\label{N2}
    \int_{0}^{\infty} D_{2} \tw(r) D_{2} w(r) dr &=
	\int_{0}^{\infty} \frac{d^{2}\tw}{dr^{2}} \frac{d^{2}w}{dr^{2}} dr
	    + \frac{3}{r}\bigl((\frac{\tw w}{r})' - \tw' w'\bigr)
		\Bigr|_{0}^{\infty} .
\end{align}
    These formulas can be derived from the definition
(\ref{Dw})
    by a straightforward integration by parts.
    It is not hard to see that in the case when
$ \tw $,
$ w $
    do not grow at infinity and vanish at zero as above expressions for
$ w^{\kappa} $,
$ w^{\lambda} $,
    the boundary terms in the RHSs disappear and only the integration
    terms survive.

    Now one can apply
(\ref{N1})
    and
(\ref{N2})
    to
(\ref{tv})
    and
(\ref{uDw})
    and calculate the norm of the (real valued) eigenfunctions
$ \tv_{l}^{\kappa} $:
\begin{align*}
    \int_{0}^{\infty} & \tv_{1}^{\kappa}(r) \tv_{1}^{\kappa}(r) dr =
	\int_{0}^{\infty} \frac{dw^{\kappa}}{dr} \frac{dw^{\kappa}}{dr} dr =\\
    &= -\kappa^{2} \int_{0}^{\infty}
	(e^{i\frac{3\pi}{4}} \exp\{e^{i\frac{3\pi}{4}}\kappa r \}
	- e^{-i\frac{3\pi}{4}} \exp\{e^{-i\frac{3\pi}{4}}\kappa r \})^{2} dr
    = \frac{\kappa}{\sqrt{2}} , \\
    \int_{0}^{\infty} & \tv_{2}^{\kappa}(r) \tv_{2}^{\kappa}(r) dr =
	\int_{0}^{\infty} \frac{d^{2}w^{\kappa}}{dr^{2}}
	    \frac{d^{2}w^{\kappa}}{dr^{2}} dr =\\
    &= \kappa^{4} \int_{0}^{\infty}
	(\exp\{e^{i\frac{3\pi}{4}}\kappa r \}
	+ \exp\{e^{-i\frac{3\pi}{4}}\kappa r \})^{2} dr
    = \frac{\kappa^{3}}{\sqrt{2}} ,
\end{align*}
    and then that of the ``continuous eigenfunctions''
$ \tu_{l}^{\lambda} $:
\begin{align*}
    &\int_{0}^{\infty} \tu_{1}^{\lambda}(r) \tu_{1}^{\mu}(r) dr =
	\int_{0}^{\infty} \frac{dw^{\lambda}}{dr}
	    \frac{dw^{\mu}}{dr} dr =\\
    &= \lambda\mu \int_{0}^{\infty}
    \bigl(\cos\lambda r+\sigma_{1}^{\lambda}(e^{-\lambda r}
	-\sin\lambda r)\bigr)
    \bigl(\cos\mu r+\sigma_{1}^{\mu}(e^{-\mu r}
	-\sin\mu r)\bigr) dr = \\
    &= \frac{\pi}{2} \lambda\mu
	(1+\sigma_{1}^{\lambda}\sigma_{1}^{\mu}) \,
	    \delta(\lambda-\mu) ,
\end{align*}
\begin{align*}
    &\int_{0}^{\infty} \tu_{2}^{\lambda}(r) \tu_{2}^{\mu}(r) dr =
	\int_{0}^{\infty} \frac{dw^{\mu}}{dr}
	    \frac{dw^{\lambda}}{dr} dr =\\
    &= \lambda^{2}\mu^{2} \int_{0}^{\infty}
    \bigl(\sin\lambda r+\sigma_{2}^{\lambda}(\cos\lambda r
	+ e^{-\lambda r})\bigr)
    \bigl(\sin\mu r+\sigma_{2}^{\mu}(\cos\mu r
	+ e^{-\mu r})\bigr) dr = \\
    &= \frac{\pi}{2} \lambda^{2}\mu^{2}
	(1+\sigma_{2}^{\lambda}\sigma_{2}^{\mu}) \,
	    \delta(\lambda-\mu) .
\end{align*}
    The above calculations allow one to write the normalized
    eigenfunctions in the following form
\begin{align*}
    v_{l}^{\kappa} &= i2^{1/4}\kappa^{1/2-l} D_{l} \bigl(
        \exp\{e^{-i\frac{3\pi}{4}} \kappa r \} -
	\exp\{e^{i\frac{3\pi}{4}} \kappa r \} \bigr), \quad \kappa > 0,\\
    u_{l}^{\kappa} &= \frac{\sqrt{2}\lambda^{-l}}{\sqrt{\pi
	(1+(\sigma_{l}^{\lambda})^{2})}} D_{l} \bigl(
	\sin\lambda r + \sigma_{l}^{\lambda} (\cos\lambda r - e^{-\lambda r})
	\bigr) ,
\end{align*}
    where the dependencies of
$ \sigma $
    on
$ \lambda $ and
$ \kappa $
    are fixed by
(\ref{S1})
    and
(\ref{S2}).
    One may notice here that the singularities present in
$ \sigma $
    for a positive
$ \kappa $
    cancel, and the overall expression for
$ u_{l}^{\kappa} $
    is well defined for any positive
$ \lambda $.

    The real valued functions
$ u_{l}^{\lambda}(r) $
    (which also depend on
$ \kappa $) and
$ v_{l}^{\kappa}(r) $
    satisfy the relations of orthonormality
\begin{align*}
    \int_{0}^{\infty} u_{l}^{\lambda_{1}}(r) u_{l}^{\lambda_{2}}(r) dr
	&= \delta(\lambda_{1}-\lambda_{2}) , \\
    \int_{0}^{\infty} u_{l}^{\lambda}(r) v_{l}^{\kappa}(r) dr
	&= 0 , \quad \kappa >0 , \\
    \int_{0}^{\infty} v_{l}^{\kappa}(r) v_{l}^{\kappa}(r) dr
	&= 0 , \quad \kappa >0 
\end{align*}
    and completeness
\begin{equation*}
    \int_{0}^{\infty} u_{l}^{\lambda}(r_{1}) u_{l}^{\lambda}(r_{2}) d\lambda
	+ v_{l}^{\kappa}(r_{1}) v_{l}^{\kappa}(r_{2}) = \delta(r_{1}-r_{2})
\end{equation*}
(the second term in the LHS is present only for the discrete spectrum, {\it i.e.} for
$ \kappa >0 $).

\section{Resolvent kernel}
\label{resolvent}
    The kernel of the resolvent is an object equivalent to the above
    spectral decomposition, but in some cases it can be 
    more appropriate to use than the latter.
    We shall look for the resolvent kernel
$ R(r,s;z) $
    of the operator
$ T_{l\kappa}^{2} $
    in form of a function which obeys the differential equation
\begin{equation}
\label{DE}
    \bigl((-\frac{d^{2}}{dr^{2}} +\frac{l(l+1)}{r^{2}})^{2} - z^{4}\bigr)
	R(r,s;z) = \delta(r-s) , \quad 0<\arg z <\frac{\pi}{2} ,
\end{equation}
    as well as the boundary conditions
(\ref{NB1}),
(\ref{NB2}),
    while being symmetric in the arguments
$ r $ and
$ s $,
    and exponentially vanishing in these arguments at infinity.
    The functions
\begin{equation*}
    D_{l} \exp\{i^{k}zr\} , \quad k=0,1,2,3
\end{equation*}
    satisfy the homogeneous equation
\begin{equation*}
    (T_{l}^{2} - z^{4})
	D_{l} \exp\{i^{k}zr\} = 0,
\end{equation*}
    where 
$ T_{l} $
    is a shorthand for the differential operation
\begin{equation*}
    T_{l} = -\frac{d^{2}}{dr^{2}} + \frac{l(l+1)}{r^{2}} .
\end{equation*}
    The exponents
$ \exp\{i^{k}zr\} $, 
$ k=0,1,2,3 $
    possess a remarkable property: in a linear combination of them,
    if the first (constant) coefficient of expansion near zero vanishes,
    then the fifth coefficient
    (the one at the power
$ r^{4} $) also vanishes (and the same, consequently, applies to all other
    coefficients at the powers multiple of 4).
    This gives an additional degree of freedom for the composition of the
    asymptotics
(\ref{NB1}) or
(\ref{NB2}).
    Specifically, this allows one to construct the functions
\begin{align*}
    h_{-} &= D_{l}(e^{-izr} + \alpha_{-}e^{izr} + \beta_{-} e^{-zr})
	\equiv \tilde{h}_{-} + \beta_{-} g_{+} , \quad g_{+} = D_{l} e^{-zr},\\
    h_{+} &= D_{l}(e^{zr} + \alpha_{+}e^{-zr} + \beta_{+} e^{izr})
	\equiv \tilde{h}_{+} + \beta_{+} g_{-} , \quad g_{-} = D_{l} e^{izr}
\end{align*}
    with the right boundary conditions, out of the three chosen exponents,
    in such a way, that functions
$ \tilde{h}_{\pm}(r) $,
$ g_{\pm}(r) $
    will satisfy the second order differential equations
\begin{equation}
\label{SOC}
    (T_{l} \pm z^{2})\tilde{h}_{\pm}(r) = 0 ,
    \quad (T_{l} \pm z^{2})g_{\pm}(r) = 0 ,
\end{equation}
    and, at the same time,
$ g_{\pm}(r) $
    will exponentially vanish at infinity.
    The coefficients
$ \alpha $ and
$ \beta $
    can be calculated directly by relating the equations
(\ref{A1}),
(\ref{A2})
    for the sum of the exponents
    to the corresponding asymptotics
(\ref{V1}) or
(\ref{V2}).
    It turns out that for
$ l=1 $:
\begin{align*}
    \alpha_{+} &= -\frac{(i+1)z+\sqrt{2}\kappa}{(i-1)z+\sqrt{2}\kappa} ,\quad
    \beta_{+} = \frac{2z}{(i-1)z+\sqrt{2}\kappa} ,\\
    \alpha_{-} &= \frac{(i+1)z-\sqrt{2}\kappa}{(i-1)z+\sqrt{2}\kappa} ,\quad
    \beta_{-} = -\frac{2iz}{(i-1)z+\sqrt{2}\kappa} ,
\end{align*}
    and for
$ l=2 $:
\begin{align*}
    \alpha_{+} &= -\frac{(i-1)z^{3}+\sqrt{2}\kappa}{(i+1)z^{3}+\sqrt{2}\kappa}
	,\quad
    \beta_{+} = -\frac{2z^{3}}{(i+1)z^{3}+\sqrt{2}\kappa} ,\\
    \alpha_{-} &= -\frac{(1-i)z^{3}+\sqrt{2}\kappa}{(i+1)z^{3}+\sqrt{2}\kappa}
	,\quad
    \beta_{-} = -\frac{2iz^{3}}{(i+1)z^{3}+\sqrt{2}\kappa} .
\end{align*}
    The functions
$ h_{\pm} $ and
$ g_{\pm} $
    introduced above allow one to contruct the resolvent kernel in the following
    way
\begin{multline}
\label{ResK}
    R(r,s;z) = \frac{1}{2z^{2}W_{-}}\bigl(
	h_{-}(r)g_{-}(s)\theta(s-r) + h_{-}(s)g_{-}(r) \theta(r-s)\bigr) -\\
    - \frac{1}{2z^{2}W_{+}}\bigl(
	h_{+}(r)g_{+}(s)\theta(s-r) + h_{+}(s)g_{+}(r) \theta(r-s)\bigr) ,
\end{multline}
    where
$ W_{\pm} $ are the wronskians
\begin{equation*}
    W_{\pm}(z) = \tilde{h}_{\pm}' g_{\pm} - \tilde{h}_{\pm} g_{\pm}' ,
\end{equation*}
    which, taking into account equations
(\ref{SOC}),
    do not depend on
$ r $
    and can be calculated at any convenient point ({\it e.g.} at infinity)
    bringing the following simple expressions
\begin{align*}
    l&=1: \quad W_{-} = -2iz^{3} , \quad W_{+} = -2z^{3} , \\
    l&=2: \quad W_{-} = -2iz^{5} , \quad W_{+} = 2z^{5} .
\end{align*}
    Symmetry of
$ R(r,s;z) $
    in the arguments
$ r $ and
$ s $
    is evident, 
    while the asymptotic conditions at zero and at infinity
    follow from the construction of the components
$ h $ and
$ g $.
    All that is left to check is the differential equation
(\ref{DE}).
    To do this, one can substitute the definitions
\begin{equation*}
    h_{+} = \tilde{h}_{+} + \beta_{+} g_{+} , \quad
    h_{-} = \tilde{h}_{-} + \beta_{-} g_{-}
\end{equation*}
    and then absorb the tilde variables
$ \tilde{h}_{\pm} $
    into terms
$ R_{\pm} $
\begin{equation*}
    R = R_{-} - R_{+} + R_{g} ,
\end{equation*}
    so that
\begin{equation*}
    R_{\pm} = \frac{1}{2z^{2}W_{\pm}}\bigl(
	\tilde{h}_{\pm}(r)g_{\pm}(s)\theta(s-r)
	    + \tilde{h}_{\pm}(s)g_{\pm}(r) \theta(r-s)\bigr),
\end{equation*}
    and
\begin{multline*}
    R_{g} = \frac{\beta_{-}}{2z^{2}W_{-}}\bigl(
	g_{+}(r)g_{-}(s)\theta(s-r) + g_{+}(s)g_{-}(r) \theta(r-s)\bigr) -\\
    - \frac{\beta_{+}}{2z^{2}W_{+}}\bigl(
	g_{-}(r)g_{+}(s)\theta(s-r) + g_{-}(s)g_{+}(r) \theta(r-s)\bigr) .
\end{multline*}
    By means of a plain substitution one can check that
\begin{equation*}
    \frac{\beta_{-}}{W_{-}} = - \frac{\beta_{+}}{W_{+}} =
	\begin{cases}
	    ((i-1)z^{3}+\sqrt{2}\kappa z^{2})^{-1} , & l=1 \\
	    ((i+1)z^{5} + \sqrt{2}\kappa z^{2})^{-1} , & l=2 ,
	\end{cases}
\end{equation*}
    separately for
$ l=1 $ and for 
$ l=2 $.
    This, together with the property
\begin{equation*}
    \theta(r-s) + \theta(s-r) = 1 ,
\end{equation*}
    yields the following expression for
$ R_{g} $
\begin{equation}
\label{Rg}
    R_{g} = \frac{\beta_{-}}{2z^{2}W_{-}} \bigl(g_{+}(r)g_{-}(s)
	+g_{+}(s)g_{-}(r)\bigr) .
\end{equation}
    This function is smooth and symmetric in its arguments
$ r $ and
$ s $,
    and evidently satisfies the homogeneous equation
\begin{equation*}
    (T_{l}^{2} - z^{4})	R_{g}(r,s;z) = 0 .
\end{equation*}

    The functions
$ R_{+} $, 
$ R_{-} $
    have the form of resolvents of second order differential operators
    constructed on the solutions
(\ref{SOC}) (although with rather strange boundary conditions).
    For that reason they satisfy the second order differential equations
\begin{equation*}
    (T_{l} \pm z^{2}) R_{\pm}(r,s;z) = \frac{1}{2z^{2}} \delta(r-s) .
\end{equation*}
    These equations allow one to write the following equality
\begin{multline*}
    (T_{l}^{2} -z^{4})(R_{-}-R_{+}) = (T_{l}+z^{2})(T_{l}-z^{2}) R_{-}
	- (T_{l}-z^{2})(T_{l}+z^{2})R_{+} = \\
    = (T_{l}+z^{2}) \frac{\delta(r-s)}{2z^{2}}
	- (T_{l}-z^{2})\frac{\delta(r-s)}{2z^{2}} = \delta(r-s) .
\end{multline*}
    This precisely means that the constructed function
$ R(r,s;z) $
    satisfies differential equation
(\ref{DE})
    with the prescribed boundary conditions, and in this way
    it represents the resolvent kernel of the self-adjoint operator
$ T_{l\kappa}^{2} $.

\subsection{Reverse operator}
    The kernel of the operator reverse to 
$ T_{l\kappa}^{2} $
    can be derived from the resolvent kernel
(\ref{ResK})
    by means of taking the limit
$ z \to 0 $.
    That calculation, however, requires an expansion to the order 5
    in parameter
$ z $
    and is rather laborious, so we shall provide a more
    primitive construction.

    Let
$ T_{l}^{-1}(r,s) $
    be a kernel of the operator reverse to
\begin{equation*}
    T_{l} = -\frac{d^{2}}{dr^{2}} + \frac{l(l+1)}{r^{2}} ,
\end{equation*}
    considered as an essentially self-adjoint
    second order differential operator.
    Then its square
\begin{equation*}
    T_{l}^{-2}(r,s) = \int_{0}^{\infty} T_{l}^{-1}(r,q) T_{l}^{-1}(q,s) dq
\end{equation*}
    does satisfy the formal equality
\begin{equation*}
    \bigl(-\frac{d^{2}}{dr^{2}} + \frac{l(l+1)}{r^{2}}\bigr)^{2}
	T_{l}^{-2}(r,s) = \delta(r-s) ,
\end{equation*}
    but does not possess the necessary asymptotic at zero.
    Then the trick is to add to the square
$ T_{l}^{-2}(r,s) $
    a term similar to
(\ref{Rg}),
    such that it, on one hand, would vanish under the action of the fourth
    order differential expression
$ T_{l}^{2} $, and on the other hand would fix the behaviour
    of
$ T_{l}^{-2}(r,s) $
    near zero.

    For
$ l=1 $
\begin{equation*}
    T_{1}^{-1}(r,s) = \frac{1}{3}\bigl(\frac{r^{2}}{s}\theta(s-r)
	+\frac{s^{2}}{r}\theta(r-s)\bigr) ,
\end{equation*}
    and after integration we obtain
\begin{equation*}
    T_{1}^{-2}(r,s) = \frac{1}{6}\bigl((r^{2}s - \frac{r^{4}}{5s})\theta(s-r)
	+(s^{2}r - \frac{s^{4}}{5r})\theta(r-s)\bigr) .
\end{equation*}
    The additional term
\begin{equation*}
    - \frac{1}{2\sqrt{2} \kappa} rs = - \frac{1}{6} \frac{3}{\sqrt{2}\kappa} rs
	\bigl(\theta(s-r)+\theta(r-s)\bigr)
\end{equation*}
    brings the expansion in
$ r $
    and 
$ s $
    near zero to the form
(\ref{V1})
\begin{multline*}
    \Theta_{1}(r,s) = T_{1}^{-2} - \frac{1}{2\sqrt{2} \kappa} rs = \\
	= \frac{1}{6}\bigl((r^{2}s - \frac{3}{\sqrt{2}\kappa}rs
	    - \frac{r^{4}}{5s})\theta(s-r)
+(s^{2}r -\frac{3}{\sqrt{2}\kappa}rs - \frac{s^{4}}{5r})\theta(r-s)\bigr) ,
\end{multline*}
    and at the same time it does not spoil the condition
\begin{equation*}
    \bigl(-\frac{d^{2}}{dr^{2}} + \frac{l(l+1)}{r^{2}}\bigr)^{2}
	\Theta_{1}(r,s) = \delta(r-s) .
\end{equation*}
    Thus we conclude that the function
$ \Theta_{1}(r,s) $
    is the kernel of the self-adjoint operator reverse to
$ T_{1\kappa}^{2} $.
    This operator is not bounded, which is manifested in the growth
    of its kernel at the infinity.

    Case
$ l=2 $
    can be treated a similar way:
\begin{equation*}
    T_{2}^{-1}(r,s) = \frac{1}{5}\bigl(\frac{r^{3}}{s^{2}}\theta(s-r)
	+\frac{s^{3}}{r^{2}}\theta(r-s)\bigr) .
\end{equation*}
    Integration gives the following expression for the kernel of the
    operator square
\begin{equation*}
    T_{2}^{-2}(r,s) = \frac{1}{10}\bigl((\frac{r^{3}}{3}
	- \frac{r^{5}}{7s^{2}})\theta(s-r)
	+(\frac{s^{3}}{3} - \frac{s^{5}}{7r^{2}})\theta(r-s)\bigr) .
\end{equation*}
    The additional term fixing the boundary conditions now is a constant
\begin{equation*}
    - \frac{1}{2 \sqrt{2} \kappa^{3}}
      = - \frac{1}{30} \frac{15}{\sqrt{2}\kappa^{3}} 
	\bigl(\theta(s-r)+\theta(r-s)\bigr) ,
\end{equation*}
    which gives the resulting expression for the kernel of the reverse
    operator as
\begin{multline*}
    \Theta_{2}(r,s) = T_{2}^{-2}(r,s) - \frac{1}{2\sqrt{2}\kappa^{3}} = \\
    = \frac{1}{10}\bigl((\frac{r^{3}}{3} - \frac{5}{\sqrt{2}\kappa^{3}}
	- \frac{r^{5}}{7s^{2}})\theta(s-r)
	+(\frac{s^{3}}{3} -\frac{5}{\sqrt{2}\kappa^{3}}
	    - \frac{s^{5}}{7r^{2}})\theta(r-s)\bigr) .
\end{multline*}

\section{Quadratic form}
\label{squadr}
    In this section we outline two expressions for the extensions
    of the quadratic form
(\ref{TVL})
    of the 3-dimensional operator.
    For that one is to substitute the transverse components
\begin{equation*}
    \vec{f}_{lm} = \sqrt{l(l+1)} \frac{u_{lm}}{r^{2}}
	\vec{\Upsilon}_{lm}(\Omega)
	+ \frac{u'_{lm}}{r}\vec{\Psi}_{lm}(\Omega)
\end{equation*}
    into the integral
\begin{equation}
\label{fint}
    \sum_{k,j} \int_{\RR^{3}\setminus B_{r}}
	|\frac{\partial f_{k}}{\partial x_{j}}|^{2} d^{3} x
\end{equation}
    over the complement of the ball
$ B_{r} $
    of the radius
$ r $
    centered at the origin. Integration by parts allows to extract the
    expressions
\begin{equation}
\label{quform}
    \int_{r}^{\infty} \ol{u_{lm}(r)}
	(-\frac{d^{2}}{dr^{2}}+\frac{l(l+1)}{r^{2}})^{2}u_{lm}(r) dr,
\end{equation}
    for the quadratic forms in question, but one should keep the terms
    with the values of the functions on the boundary of
$ B_{r} $.
    By extending the set of the parametrizing functions
$ u_{lm} $
    to the domain of
$ T_{l\kappa}^{2} $
    with the boundary conditions
(\ref{NB1}),
(\ref{NB2}),
    the expressions
(\ref{quform})
    are transferred into the extended quadratic forms of these operators.
    In this process the second (for $ l=1 $) and the third (for $ l=2 $)
    derivative in the boundary terms is to be expressed via the first
    derivative
    or the value of the function, in correspondence with
(\ref{NB1}) or
(\ref{NB2}).
    Further the boundary terms can be expressed via
$ \vec{f}(\vec{x}) $
    and gathered into
\begin{multline*}
    f_{lm}^{2}(r) = \sum_{k,k'}\iint_{\Sph^{2}} f_{k}(r,\Omega) \bigl(
\ol{\Upsilon^{k}(\Omega)} \Upsilon^{k'}(\Omega')
+\ol{\Psi^{k}(\Omega)} \Psi^{k'}(\Omega')
    \bigr) \times\\
\times \ol{f_{k'}(r,\Omega')}
    r^{2} d\Omega d\Omega'
	= |u'_{lm}(r)|^{2} + \frac{l(l+1)}{r^{2}} |u_{lm}(r)|^{2} .
\end{multline*}
    Finally they enter the extended quadratic form in the following way
\begin{multline*}
    Q_{\kappa}(f) = \lim_{r\to 0}\Bigl(
    \int_{\RR^{3}\setminus B_{r}}
	|\frac{\partial f_{k}}{\partial x_{j}}|^{2} d^{3} x -\\
- \bigl(\frac{22\sqrt{2}}{9}\kappa_{1m}+\frac{5}{3r}\bigr)f_{1m}^{2}(r) -
\bigl(\frac{80\sqrt{2}}{750}\kappa_{1m}^{3}+\frac{4}{r}\bigr)f_{2m}^{2}(r)
    \Bigr) ,
\end{multline*}
    where
$ k $,
$ j $ and
$ m $
    are summed up and
$ \vec{f}(\vec{x}) $
    is transverse as in
(\ref{transc}).
    One can notice, that for regular
$ \vec{f}(\vec{x}) $
    the boundary terms vanish, the limit
$ r\to 0 $
    of the integral
(\ref{fint}) is finite
    and the quadratic form
$ Q_{\kappa}(f) $
    coincides with
(\ref{TVL}).
    The expression for
$ Q_{\kappa}(f) $
    is nontrivial when some of the components of
$ \vec{f}(\vec{x}) $
    with angular momentum
$ l=1 $ or
$ l=2 $
    diverge at the origin as
$ r^{-1} $
    or
$ r^{-2} $,
    correspondingly.
    In this case the limit of the integral
(\ref{fint})
    over the complement of the ball
$ B_{r} $
    diverges as
$ r\to 0 $
    and the boundary terms are arranged in the way to produce a finite
$ \kappa $-dependent 
    value for
$ Q_{\kappa}(f) $.
    In the sperical symmetric case involving only the
$ l=1 $
    components ({\it i.e.} for $\kappa_{1m} = \kappa $, 
$ \kappa_{2m} = -\infty$),
    the expression for
$ Q_{\kappa}(f) $
    can be rewritten in a simpler form similar to
(\ref{Qfs})
\begin{equation*}
    Q_{\kappa}(f) = \lim_{r\to 0}\Bigl(
    \int_{\RR^{3}\setminus B_{r}}
	|\frac{\partial f_{k}}{\partial x_{j}}|^{2} d^{3} x -
    (\frac{5}{3r}+ \frac{22\sqrt{2}}{9}\kappa) \int_{\partial B_{r}}
	|\vec{f}(\vec{x})|^{2} d^{2} s \Bigr),
\end{equation*}
    where the function
$ \vec{f}(\vec{x}) $
    can be singular as
$ r^{-1} $
    at the origin, but is still square integrable over
$ \RR^{3} $.

\section{Conclusion}
    We have discussed self-adjoint extensions of the fourth-order differential
    operators
$ T_{l0}^{2} $,
$ l=1,2 $,
    acting on the parametrizing functions
    in the quadratic form of the transverse vector Laplace operator.
    These operators allow to contruct extentions of the original quadratic
    form on functions in the 3-dimensional space.
    One has
$ 3+5=8 $
    dimensional extension parameters
$ \kappa_{lm} $,
    which in a special case can be equal, thus setting up an additional symmetry
    of the model. In particular, at
$ \kappa_{lm} = \kappa > 0 $
    the extended quadratic form
$ Q_{\kappa}(\vec{f}) $
    possesses stable (with respect to variation) states which look as
\begin{equation*}
    \sum_{|m|\leq l\leq 2} A_{lm} \bigl(
	\frac{v_{l}^{\kappa}(r)}{r^{2}}\vec{\Upsilon}_{lm}(\Omega)
	+\frac{v_{l}^{\kappa}{}'(r)}{r}\vec{\Psi}_{lm}(\Omega) \bigr) ,
\end{equation*}
    wherein
$ A_{lm} $
    can be either scalars or elements of the representation 
    of some algebra of internal symmetry.

    It is worthwhile to note that the introduction of nontrivial dimensional
    extension parameters breaks the dilation (scale) homogeneity 
    of the quadratic form
(\ref{TVL}).
    That is, the presented extension
    can be treated as a way to introduce
    a dimensional parameter in the model of \emph{classical} mechanics
    with the potential energy of the form
(\ref{TVL}).

\section*{Acknowledgments}
    The author is grateful to P.~A.~Bolokhov and S.~Derkachov
    for discussions. The work is partially supported by
RFBR grants 14-01-00341, 12-01-00207
    and the programme ``Mathematical problems of nonlinear dynamics'' of RAS.


\begin{thebibliography}{0}

\bibitem{FS} K.~Friedrichs, ``Spektraltheorie halbbeschr\"ankter Operatoren,''
    Math. Ann. {\bf 109}, 1934, 465--487;\\
    M.~Stone, in \emph{Linear Transformations in Hilbert spaces and their
    Applications in Analysis}, Amer. Math. Soc. Colloquium Publication {\bf 15},
    Providence, R.I., 1932;\\
    or see theorem X.23 in \cite{RS}.

\bibitem{RS} M.~Reed, B.~Simon, \emph{Methods of Modern Mathematical Physics.
II: Fourier Analysis, Self-adjointness}, Academic Press, 1975.

\bibitem{Krein}
    M.~G.~Krein, ``The theory of self-adjoint extensions of semi-bounded
Hermitian transformations and its applications.'',
    Rec. Math. (Mat. Sbornik) N.S., {\bf 20} (62), 1947, 431--495;
    Rec. Math. (Mat. Sbornik) N.S., {\bf 21} (64), 1947, 365--404.

\bibitem{BF}
    F.~A.~Berezin, L.~D.~Faddeev,
  ``A Remark on Schrodinger's equation with a singular potential,''
  Sov.\ Math.\ Dokl.\  {\bf 2} (1961) 372
  [Dokl.\ Akad.\ Nauk Ser.\ Fiz.\  {\bf 137} (1961) 1011].

\bibitem{AK}
    S.~Albeverio, P.~Kurasov, \emph{Singular Perturbation of Differential
    Operators. Solvable Schr\"odinger type Operators},
    Cambridge University Press, 2000.

\bibitem{VSH}
    see e.g. E.~L.~Hill, ``The Theory of Vector Spherical Harmonics'',
    Am. J. Phys. {\bf 22} (1954) 211.

\bibitem{Richt}
    R.~D.~Richtmyer, \emph{Principles of Advanced Mathematical Physics, vol.1},
    Springer-Verlag, New York Heildelberg Berlin, 1978.
    

\end{thebibliography}
\end{document}